\documentclass[smallcondensed]{svjour3}

\usepackage[american]{babel}
\usepackage[utf8]{inputenc}
\usepackage[T1]{fontenc}
\usepackage[binary-units=true]{siunitx}
\usepackage{booktabs}
\usepackage{xspace}
\usepackage{tikz}
\usepackage{microtype}
\usepackage{todonotes}
\usepackage{csquotes}
\usepackage{amssymb}
\usepackage{mathtools}
\usepackage{mathptmx}

\smartqed 

\newcommand{\graphset}{}

\newcommand{\Graph}{\graphset{G}}
\newcommand{\Arcs}{\graphset{A}}

\newcommand{\Nodes}{\graphset{N}}
\newcommand{\Vertices}{\Nodes}

\newcommand{\ematpH}{Q}

\newcommand{\cellind}{\beta}

\newcommand{\othcellind}{\sigma}
\newcommand{\numcells}{\kappa}
\newcommand{\cellset}{V}

\newcommand{\node}{n}
\newcommand{\otherNode}{m}

\newcommand{\arc}{a}
\newcommand{\edge}{\arc}
\newcommand{\otheredge}{b}
\newcommand{\pipind}{\alpha}
\newcommand{\othpipind}{\mu}

\newcommand{\Outedges}{\delta^{\text{out}}}

\newcommand{\Inedges}{\delta^{\text{in}}}

\newcommand{\Vdh}{\Nodes}

\newcommand{\Vff}{\Nodes_{\text{ff}}}
\newcommand{\Vbf}{\Nodes_{\text{bf}}}
\newcommand{\Aff}{\Arcs_{\text{ff}}}
\newcommand{\Abf}{\Arcs_{\text{bf}}}
\newcommand{\Ac}{\Arcs_{\text{c}}}
\newcommand{\ad}{\arc_{\text{d}}}
\newcommand{\inflowArcs}{\mathcal{I}}
\newcommand{\outflowArcs}{\mathcal{O}}
\newcommand{\pipes}{\Arcs_{\text{p}}}

\makeatletter
\newcommand{\fcdot}{\,\cdot\,}
\newcommand{\fcarg}[1]{\def\fc@rg{#1}\ifx\fc@rg\empty\fcdot\else\fc@rg\fi}
\makeatother

\newcommand{\defset}[3][\defsep]{\set{#2#1#3}}

\newcommand{\set}[1]{\{#1\}}

\newcommand{\feedin}{P_{\text{in}}}
\newcommand{\edens}{e}
\newcommand{\volflo}{\hat q}

\newcommand{\temperature}{T}
\newcommand{\temp}{\temperature}
\newcommand{\massflow}{q}
\newcommand{\mflow}{\massflow}
\newcommand{\pressure}{p}
\newcommand{\press}{\pressure}
\newcommand{\heatcap}{c_{\text{p}}}

\newcommand{\Plim}{\bar \Power}
\newcommand{\tempbf}{\temp^{\text{bf}}}

\newcommand{\maxtempff}{\temp^{\text{ff}}_+}
\newcommand{\mintempff}{\temp^{\text{ff}}_-}

\newcommand{\maxpresff}{\press^{\text{ff}}_+}
\newcommand{\minpresff}{\press^{\text{ff}}_-}

\newcommand{\maxpresbf}{\press^{\text{bf}}_+}
\newcommand{\minpresbf}{\press^{\text{bf}}_-}

\newcommand{\tempdiffhs}{\Delta \temp}

\newcommand{\timeint}{[t_0,t_{\text{end}}]}

\newcommand{\length}{\ell}
\newcommand{\velocity}{v}
\newcommand{\vel}{\velocity}

\newcommand{\Power}{P}

\newcommand{\abbr}[1][abbrev]{#1.\xspace}

\newcommand{\cf}{\abbr[cf]}
\newcommand{\eg}{\abbr[e.g]}

\newcommand{\ie}{\abbr[i.e]}

\journalname{Differential-Algebraic Equations Forum}

\begin{document}

\title{Port-Hamiltonian modeling of district heating
  networks}

\author{Sarah-Alexa Hauschild \and Nicole Marheineke
  \and Volker Mehrmann \and Jan Mohring \and Arbi Moses Badlyan \and
  Markus Rein \and Martin Schmidt}

\authorrunning{Hauschild et al.}

\institute{S. Hauschild \and N. Marheineke \and M. Schmidt \at
  Trier University, Department of Mathematics, Universitätsring 15,
  D-54296 Trier, Germany \\
  \email{\{hauschild, marheineke, martin.schmidt\}@uni-trier.de}
  \and
  V. Mehrmann \and A. Moses Badlyan \at
  TU Berlin, Institut f\"ur Mathematik, MA 4-5, Straße des 17. Juni 136, D-10587 Berlin, Germany\\
  \email{\{mehrmann, badlyan\}@math.tu-berlin.de}
  \and
  J. Mohring \and M. Rein \at
  Fraunhofer-ITWM, Fraunhofer-Platz 1, D-67663 Kaiserslautern, Germany\\
  \email{\{jan.mohring, markus.rein\}@itwm.fraunhofer.de}
}

\date{Received: date / Accepted: date}

\maketitle

\begin{abstract}
  This paper provides a first contribution to port-Hamiltonian modeling of district heating networks. By introducing a model hierarchy of flow equations on the network, this
work aims at a thermodynamically consistent port-Hamiltonian embedding
of the partial differential-algebraic systems. We show that
a spatially discretized network model describing the advection of the
internal energy density with respect to an underlying incompressible
stationary Euler-type hydrodynamics can be considered as a
parameter-dependent finite-dimensional port-Hamiltonian
system. Moreover, we present an infinite-dimensional port-Hamiltonian
formulation for a compressible instationary thermodynamic fluid flow in a pipe. Based on these first promising results, we raise open questions and point out research
perspectives concerning structure-preserving discretization, model
reduction, and optimization.




\end{abstract}

\section{Introduction}
\label{sec:introduction}

A very important part of a successful energy transition is an
increasing supply of renewable energies.
However, the power supply through such energies is highly volatile.
That is why a balancing of this volatility and more energy efficiency
is needed. An important player in this context are district heating networks.
They show a high potential to balance the fluctuating supply of
renewable energies due to their ability to absorb more or less excess power while keeping the heat supply unchanged. A long-term objective is to strongly increase energy
efficiency through the intelligent control of district heating
networks.
The basis for achieving this goal is the dynamic modeling of the
district heating network itself, which is not available in the optimization tools
currently used in industry. Such a dynamic modeling would allow for
optimization of the fluctuating operating resources, \eg, waste
incineration, electric power, or gas. However, as power and heating networks act
on different time scales and since their descriptions lead to mathematical
problems of high spatial dimension, their coupling for a dynamic
simulation that is efficiently realizable involves various
mathematical challenges. One possible remedy is a port-Hamiltonian modeling framework: Such an energy-based formulation
brings the different scales on a single level, the port-Hamiltonian
character is inherited during the coupling of individual systems, and
in a port-Hamiltonian system the physical principles (stability,
passivity, conservation of energy and momentum) are ideally encoded in
the algebraic and geometric structures of the model. Deriving model
hierarchies by using adequate Galerkin projection-based techniques for
structure-preserving discretization as well as model reduction, and
combining them with efficient adaptive optimization strategies opens up a new promising approach to complex application issues.

Against the background of this vision, this paper provides a first
contribution to port-Hamiltonian modeling of district heating
networks, illustrating the potential for optimization in a case study,
and raising open research questions and challenges. Port-Hamiltonian
(pH) systems have been elaborately studied in literature lately; see,
\eg, \cite{BeaMXZ18,MehMW18,SchJ14,SchM18} and the references therein.
The standard form appears as
\begin{subequations}\label{eq:pH}
  \begin{align}
    \frac{\mathrm{d}{z}}{\mathrm{d}t}
    = (J-R) \nabla_z \mathcal{H}(z) + (B-P) u,
    \quad
    y = (B+P)^T \nabla_z \mathcal{H}(z) + (S+N) u
  \end{align}
  with
  \begin{align}
    W = W^T\geq 0,
    \quad W =\left[
    \begin{array}{cc}
      R & P \\ P^T & S
    \end{array}\right].
  \end{align}
\end{subequations}
The Hamiltonian $\mathcal{H}$ is an energy storage function, $J =
-J^T$ is the structure matrix describing energy flux among energy
storage elements, $R = R^T$ is the dissipation matrix, $B\pm P$ are
port matrices for energy in- and output, and $S=S^T$, $N=-N^T$ are
matrices associated with the direct feed-through from input $u$ to
output $y$. The system satisfies a dissipation inequality, which is an
immediate consequence of the positive (semi-)definiteness of the
passivity matrix~$W$ and also holds even when the coefficient matrices
depend on the state $z$ or explicitly on time $t$, or when they are
defined as linear operators on infinite-dimensional spaces. Including
time-varying state constraints yields a port-Hamiltonian descriptor
system of differential-algebraic equations \cite{BeaMXZ18,MehMW18,Sch13}. Port-Hamiltonian partial
differential equations on networks (port-Hamiltonian PDAE) are topic
in, \eg, \cite{a28:egger:2018} for linear damped wave equations or
in~\cite{c47:liljegren-sailer:2019} for nonlinear isothermal Euler
equations. The adequate handling of thermal effects is a novelty of
this work. Extending the work of \cite{BadMBM18,BadZ18}, we make use of a thermodynamically consistent generalization of
the port-Hamiltonian framework in which the Hamiltonian is combined with
an entropy function. The resulting dynamic system consists of a
(reversible) Hamiltonian system and a generalized (dissipative)
gradient system. Degeneracy conditions ensure that the flows of the
two parts do not overlap. Respective pH-models in operator form can be
found, \eg, for the Vlasov--Maxwell system in plasma physics in
\cite{KraH17,KraKMS17}, for the Navier--Stokes equations for reactive flows in \cite{AltS17} or for finite strain thermoelastodynamics in \cite{BetS19}.

The paper is structured as follows. Starting with the description of a
district heating network as a connected and directed graph in
Sect.~\ref{sec:network-modeling}, we present models associated to the
arcs for the pipelines, consumers, and the depot of the network
operator that are coupled with respect to conservation of mass and energy
as well as continuity of pressure at the network's nodes. We especially
introduce a hierarchy of pipe models ranging from the compressible
instationary Navier--Stokes equations for a thermodynamic fluid flow to
an advection equation for the internal energy density coupled with
incompressible stationary Euler-like equations for the
hydrodynamics. Focusing on the latter, we show that the associated
spatially discretized network model can be embedded into a family of
parameter-dependent standard port-Hamiltonian systems in
Sect.~\ref{sec:ph_semidiscrete} and numerically explore the network's behavior in Sect.~\ref{sec:case_study}. In a study on operating the
heating network with respect to the avoidance of power peaks in the
feed-in, we particularly reveal the potential for optimization. In view
of the other pipe models, a generalization of the port-Hamiltonian
framework to cover the dissipative thermal effects is necessary. In
Sect.~\ref{sec:pH-modeling} we develop an infinite-dimensional
thermodynamically consistent port-Hamiltonian formulation for the one-dimensional partial differential equations of a compressible instationary turbulent pipe flow. From this, we raise open research questions and perspectives concerning structure-preserving discretization, model reduction, and optimization in Sect.~\ref{sec:concl-persp}.


\section{Network modeling}
\label{sec:network-modeling}

\begin{figure}[tb]
  \centering
  \begin{tikzpicture}
  \def \radius {11pt}
  \def \locations {3/0/1, 4.5/1/2, 4.5/-1/3, 6/0/4, 6/-2/5}
  \def \distance {1}
  \def \lineThickness {very thick}

  \node[circle, draw, minimum width=\radius](s) at (0,0){ };
  \foreach \x/\y/\name in \locations
  \node[circle, draw, minimum width=\radius](\name) at (\x,\y){ };

  \node[circle, draw, minimum width=\radius](t)
  at (0,0 - \distance){};
  \foreach \x/\y/\name in \locations
  \node[circle, draw, minimum width=\radius](R\name) at (\x,\y -
  \distance){ };

  \foreach \source/\dest in {R5/R3, R4/R3, R4/R2, R3/R1, R2/R1, R1/t}
  \draw [->, dashed, \lineThickness] (\source) -- (\dest);

  \foreach \source/\dest in {s/1, 1/2, 1/3, 2/4, 3/4, 3/5}
  \draw [->, \lineThickness] (\source) -- (\dest);

  \foreach \n in {1,2,4,5}
  \draw [->, blue, dotted, \lineThickness] (\n) -- (R\n);

  \draw [->, red, dashdotted, \lineThickness] (t) -- (s);
\end{tikzpicture}
%
  \caption{A schematic district heating network: Foreflow arcs are
    plotted in solid black, backflow arcs in dashed black, consumers (households) in
    dotted blue, and the depot in dash-dotted red.}
  \label{fig:sample-network}
\end{figure}

The district heating network is modeled by a connected and directed
graph $\Graph = (\Vertices, \Arcs)$ with node set~$N$ and arc
set~$A$. This graph consists of (i) a
foreflow part, which provides the consumers with hot water; (ii)
consumers, that obtain power via heat exchangers; (iii) a backflow part,
which transports the cooled water back to the depot; and (iv) the
depot, where the heating of the cooled water takes place; see
Fig.~\ref{fig:sample-network} for a schematic illustration.
The nodes $\Vdh = \Vff \cup \Vbf$ are
the disjoint union of nodes~$\Vff$ of the foreflow part and
nodes~$\Vbf$ of the backflow part of the network. The arcs $\Arcs =
\Aff \cup \Abf \cup \Ac \cup \set{\ad}$ are divided into foreflow arcs
$\Aff$, backflow arcs $\Abf$, consumer arcs $\Ac$, and the depot arc
$\ad$ of the district heating network operator. The set of pipelines
is thus given by $\pipes = \Aff \cup \Abf$.

In the following we introduce a model hierarchy for the flow in a
single pipe (\cf Fig.~\ref{fig:pipe_models}) and afterward discuss the nodal coupling conditions for the
network. Models for consumers (households) and the depot yield the
closure conditions for the modeling of the network.

\subsection{Model hierarchy for pipe flow}
\label{subsec:pipe-modeling}

Let $a \in \pipes$ be a pipe. Starting point for the modeling of the flow in a
pipe are the cross-sectionally averaged one-dimensional instationary
compressible Navier--Stokes equations for a thermodynamic fluid flow \cite{Schlichting06}.
We assume that the pipe is cylindrically shaped,
that it has constant circular cross-sections, and that the flow
quantities are only varying along the cylinder axis. Consider $(x,t)\in
(0,\ell)\times (t_0,t_\mathrm{end}] \subseteq \mathbb{R}^2$
with pipe length $\ell$ as well as start and end time $t_0$,
$t_\text{end}>0$. Mass density, velocity, and internal energy density,
\ie, $\rho, v, e: (0,\ell) \times (t_0,t_\mathrm{end}] \rightarrow
\mathbb{R}$, are then described by the balance equations
\begin{equation}
  \label{eq:distr-heat-pipe}
  \begin{split}
    0&=\partial_t\rho+\partial_x(\rho v),\\
    0&=\partial_t(\rho v)+\partial_x(\rho v^2)+\partial_x p+\frac{\lambda}{2d}\rho|v|v+\rho g\partial_x h,\\
    0&=\partial_t e+\partial_x(ev)+p\partial_xv-\frac{\lambda}{2d}\rho|v|v^2+\frac{4k_\mathrm{w}}{d}(T-\vartheta).
  \end{split}
\end{equation}
Pressure and temperature, \ie, $p, T : (0,\ell) \times
(t_0,t_\mathrm{end}] \rightarrow \mathbb{R}$, are determined by
respective state equations. In the momentum balance the frictional
forces with friction factor $\lambda$ and pipe diameter $d$ come from
the three-dimensional surface conditions on the pipe walls, the outer
forces arise from gravity with gravitational acceleration $g$ and pipe
level $h$ (with constant pipe slope $\partial_x h$). The energy exchange with the outer surrounding is
modeled in terms of the pipe's heat transmission coefficient $k_\mathrm{w}$ and
the outer ground temperature
$\vartheta$. System~\eqref{eq:distr-heat-pipe} are 
(Euler-like) non-linear hyperbolic partial differential equations of first order for a turbulent pipe flow.

The hot water in the pipe is under such a high pressure that it does
not turn into steam. Thus, the transition to the incompressible limit
of \eqref{eq:distr-heat-pipe} makes sense, yielding the following
partial differential-algebraic system for velocity $v$ and internal
energy density $e$, where the pressure $p$ acts as a Lagrange
multiplier to the incompressibility constraint:
\begin{equation}
  \label{eq:distr-heat-pipe-incomp}
  \begin{split}
    0&=\partial_x\vel,\\
    0&=\partial_t\vel+\frac{1}{\rho}\partial_x\press+\frac{\lambda}{2d}|\vel|\vel+g\partial_x
    h,\\
    0&=\partial_t e+\vel\partial_x e-\frac{\lambda}{2d}\rho|v|v^2+\frac{4k_\mathrm{w}}{d}(\temp-\vartheta).
  \end{split}
\end{equation}
The system is supplemented with state equations for density~$\rho$ and
temperature~$T$. Note that the energy term due to friction is
negligibly small in this case and can be omitted.

Since the hydrodynamic and thermal effects act on different time
scales, System~\eqref{eq:distr-heat-pipe-incomp} may be simplified
even further by setting $\partial_tv=0$, \ie,
\begin{equation}\label{eq:distr-heat-pipe-smallacc}
  \begin{split}
    0&=\partial_x\vel,\\
    0&=\partial_x\press+ \frac{\lambda}{2d}\rho |\vel|\vel+ \rho g\partial_x h,\\
    0&=\partial_t e+\vel\partial_x e+\frac{4k_\mathrm{w}}{d}(\temp-\vartheta),
  \end{split}
\end{equation}
again supplemented with state equations for $\rho$, $T$.
System~\eqref{eq:distr-heat-pipe-smallacc} describes the heat
transport in the pipe where flow velocity and pressure act as
Lagrange multipliers to the stationary hydrodynamic
equations. However, the flow field is not stationary at all because of
the time-dependent closure (boundary) conditions (at households and
the depot). In the presented model hierarchy one might even go a step
further and ignore the term concerning the heat transition with the
outer surrounding of the pipe, \ie, $4 k_\mathrm{w}(T-\vartheta)/d=0$, when
studying the overall network behavior caused by different operation of
the depot; see Sect.~\ref{sec:ph_semidiscrete} and
Sect.~\ref{sec:case_study}.

\tikzstyle{block} = [rectangle, draw, 
text width=30em, text centered, rounded corners, minimum height=3em]
\tikzstyle{line} = [draw, ->]

\begin{figure}[t]
\centering
  \begin{tikzpicture}[node distance = 1.75cm, auto]
    \node [block] (Comp) {compressible instationary thermodynamic turbulent flow  \eqref{eq:distr-heat-pipe}};
    \node [block, below of = Comp] (Incomp) {incompressible instationary thermodynamic turbulent flow \eqref{eq:distr-heat-pipe-incomp}};
    \node [block, below of = Incomp] (Stat) {energy advection with outer cooling\\ w.r.t. incompressible stationary hydrodynamic equations (\ref{eq:distr-heat-pipe-smallacc})};
     \node [block, below of = Stat] (WC) {energy advection without outer cooling \\w.r.t. incompressible stationary hydrodynamic equations (\ref{eq:e_net})};
    \path [line] (Comp) -- node{$\partial_xv=0$} (Incomp);
    \path [line] (Incomp) -- node{$\partial_tv=0$, $\frac{\lambda}{2d}\rho|v|v^2$ small}(Stat);
    \path [line] (Stat) -- node{$\frac{4k_W}{d}(T-\vartheta)=0$}(WC);
  \end{tikzpicture}
\caption{Hierarchy of pipe flow models \label{fig:pipe_models}}
\end{figure}
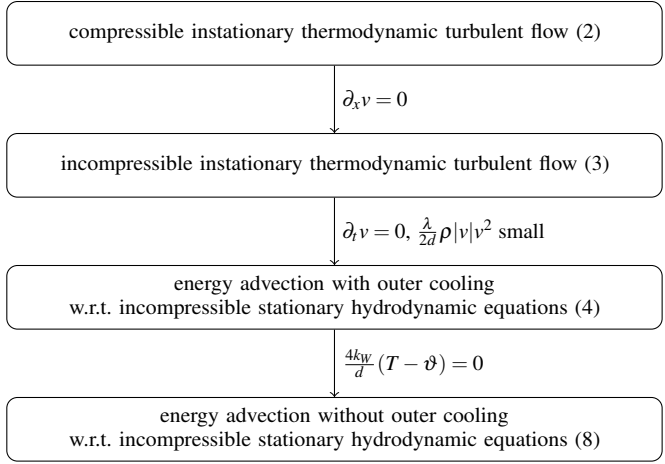

\paragraph{State equations and material models }
In the pressure and temperature regime being relevant for operating
district heating networks, we model the material properties of water by polynomials depending exclusively on the internal energy
density, and not on the pressure. The relations for temperature $T$,
mass density $\rho$, and kinematic viscosity $\bar\nu$ summarized in
Table~\ref{tab:wasser} are based on a fitting of data taken from the
NIST Chemistry WebBook~\cite{nist16}. The relative error of the
approximation is of order $O(10^{-3})$, which is slightly higher than
the error $O(\num{5e-4})$ we observe due to neglecting
the pressure dependence. The quadratic state equation for the
temperature allows a simple conversion between $e$ and $T$, which is
necessary since closure conditions (households, depot) are usually
stated in terms of $T$;
cf. Sect.~\ref{sec:network-modeling:househelds-depot}. Obviously,
$e_\star(T_\star)= 0.5 \, T_2^{-1}(-T_1 + (T_1^2 - 4 T_2( T_0 -
T_\star))^{1/2})$ holds for $T_\star(e_\star)=\sum_{i=0}^2 T_i
e_\star^i$, $e_\star \geq 0$.

\begin{table}[t]
  \centering
  \caption{Material properties of water as functions of the internal
    energy density $z(e)=z_0\,z_\star(e/e_0)$, $z
    \in\{T,\rho,\bar\nu\}$, where $z_\star$ denotes the dimensionless
    quantity scaled with the reference value $z_0$; in particular $e_\star=e/e_0$ with
    $e_0= 10^9\,\si{\joule\per\cubic\meter}$.  The stated relative errors of the underlying
    polynomial approximation hold in the regime $e \in [0.2, 0.5]$\,\si{\giga\joule\per\cubic\meter}
    and $p\in [5, 25]$\,\si{\bar}, implying $T\in [50, 130]$\,\si{\degreeCelsius}.}
  \label{tab:wasser}
  \begin{tabular}{lll}
    \toprule
    Reference
    & Material model
    &  Rel. error\\
    \midrule
    $T_0=1$\,\si{\degreeCelsius}
    & $T_\star(e_\star)=59.2453 \, e_\star^2 + 220.536 \, e_\star +
      1.93729$
    & \num{1.2e-3} \\
    $\rho_0=10^3$\,\si{\kilogram\per\cubic\meter}\hspace*{-0.1cm}
    & $\rho_\star(e_\star)=-0.208084 \, e_\star^2 -0.025576\, e_\star +
      1.00280$
    &  \num{6.0e-4}\\
    $\bar\nu_0=10^{-6}$\,\si{\square\meter\per\second}\hspace*{-0.2cm}
    &  $\bar\nu_\star(e_\star)=11.9285\, e_\star^4 -22.8079\, e_\star^3
      +17.6559\, e_\star^2 -7.00355\, e_\star +1.42624$\hspace*{-0.1cm}
    & \num{9.9e-4}\\
    \bottomrule
  \end{tabular}
\end{table}

\begin{remark}\label{rem:state}
  Alternatively to the specific data-driven approach, the state equations
  can be certainly also deduced more rigorously from thermodynamic
  laws. A thermodynamic fluid flow described by
  \eqref{eq:distr-heat-pipe} satisfies the entropy balance for
  $s:(0,\ell) \times (t_0,t_\mathrm{end}] \rightarrow \mathbb{R}$, \ie,
  \begin{align*}
    0=\partial_t s+\partial_x(sv)-\frac{\lambda}{2d}\frac{1}{\temp} \rho
    |\vel|\vel^2+\frac{4k_\mathrm{w}}{d}\frac{1}{\temp}(\temp-\vartheta).
  \end{align*}
  Considering the entropy as a function of mass density and internal
  energy density, $s=s(\rho, e)$, yields the Gibbs identities which can
  be used as state equations for pressure~$p$ and temperature~$T$, \ie,
  \begin{align*}
    \partial_\rho s=-(\rho T)^{-1} (e+p-Ts), \quad \partial_e s =T^{-1}.
  \end{align*}
\end{remark}

\paragraph{Pipe-related models }
The pipe flow is mainly driven in a turbulent regime, \ie, with
Reynolds number $\mathrm{Re}>10^3$. Thus, the pipe friction factor $\lambda$ can be
described by the Colebrook--White equation in terms of the Reynolds
number $\mathrm{Re}$ and the ratio of pipe roughness and diameter
$k_\mathrm{r}/d$,
\begin{equation*}
  \frac{1}{\sqrt{\lambda}} (v,e)
  = -2\, \log_{10} \left(
    \frac{2.52}{\mathrm{Re}(v,e)\,
      \sqrt{\lambda}(v,e)}
    + \frac{1}{3.71}\frac{k_\mathrm{r}}{d} \right),
  \quad
  \mathrm{Re}(v,e) = \frac{|v|\, d}{\bar \nu(e)}.
\end{equation*}
The model is used for technically rough pipes.
Its limit behavior corresponds to the relation by Prandtl and Karman
for a hydraulically smooth pipe, \ie, $1/{\sqrt{\lambda}} = 2
\log_{10}(\mathrm{Re}\sqrt{\lambda})-0.8$ for $k_\mathrm{r}/d\rightarrow 0$,
and to the relation by Prandtl, Karman, and Nikuradse for a completely
rough pipe, \ie, $1/{\sqrt{\lambda}}=1.14-2 \log_{10}(k_\mathrm{r}/d)$ for
$\mathrm{Re}\rightarrow \infty$, \cite{Shashi15}.
The underlying root finding problem for~$\lambda$ can be
solved using the Lambert W-function; see~\cite{clamond2009efficient}.
However, in view of the computational effort it can also be reasonable
to consider a fixed constant Reynolds number for the pipe as further simplification.

The pipe quantities -- length~$\ell$, diameter~$d$, slope~$\partial_x
h$, roughness~$k_\mathrm{r}$, and heat transmission coefficient~$k_\mathrm{w}$ -- are assumed to be constant in the pipe model. Moreover, note that in this work we also consider the outer ground temperature $\vartheta$ as constant, which will play a role for our port-Hamiltonian formulation of \eqref{eq:distr-heat-pipe} in Sect.~\ref{sec:pH-modeling} .

\subsection{Nodal coupling conditions}
\label{sec:network-modeling:nodal-coupl-cond}

For the network modeling it is convenient to use the following
standard notation. Quantities related to an arc
$\arc = (\otherNode, \node) \in \Arcs$, $\otherNode, \node \in \Vertices$,
are marked with the subscript~$\arc$, quantities associated to a node
$\node \in \Vertices$ with the subscript~$\node$. For a node $\node
\in \Vertices$, let $\Inedges_\node$, $\Outedges_\node$ be the sets of
all topological ingoing and outgoing arcs, \ie,
\begin{equation*}
  \Inedges_\node
  = \defset{\arc \in \Arcs}{\exists \otherNode \text{ with
    }\arc = (\otherNode, \node)},
  \quad
  \Outedges_\node
  = \defset{\arc \in \Arcs}{\exists
    \otherNode \text{ with } \arc = (\node, \otherNode)},
\end{equation*}
and let $\inflowArcs_\node(t)$, $\outflowArcs_\node(t)$, $t \in \timeint$,
be the sets of all flow-specific ingoing and outgoing arcs,
\begin{align*}
  \inflowArcs_n(t)
  &=\defset{\arc \in \Inedges_\node}{q_a(\ell_a,t)\geq 0}\cup
    \defset{\arc \in \Outedges_\node}{q_a(0,t)\leq 0},\\
  \outflowArcs_n(t)
  &=\defset{\arc \in \Inedges_\node}{q_a(\ell_a,t)< 0}\cup
    \defset{\arc \in \Outedges_\node}{q_a(0,t)> 0};
\end{align*}
see, \eg,
\cite{Geissler_et_al:2015,Geissler_et_al:2018,Hante_Schmidt:2019}
where a similar notation is used in the context of gas networks.
Note that the sets~$\inflowArcs_n(t)$, $\outflowArcs_n(t)$ depend on
the flow~$q_a$, $a \in A$, in the network, which is not known a
priori.

The coupling conditions we require for the network ensure the conservation of mass
and energy as well as the continuity of pressure at every node $\node\in
\Vertices$ and for all time $t\in \timeint$, i.e.,
\begin{subequations}
  \label{eq:nodal_cond}
  \begin{align}
    \label{eq:mass-balance}
    \sum_{\arc \in \Inedges_\node} \mflow_\arc(\length_\arc, t) &=
                                                                  \sum_{\arc \in \Outedges_\node} \mflow_\arc(0, t),\\
    \label{eq:energy-balance}
    \sum_{\arc \in \Inedges_\node} \volflo_\arc(\length_\arc, t) \edens_\arc(\length_\arc, t) &=
                                                                                                \sum_{\arc \in \Outedges_\node} \volflo_\arc(0, t) \edens_\arc(0, t),&& \,
                                                                                                                                                                        \edens_\arc(0, t) = \edens_\node(t), \,\,\,\, \arc \in \outflowArcs_\node(t),\\
    \label{eq:pressure-continuity}
    \press_\arc(\length_\arc, t)
                                                                &=  \press_\node(t),
                                                                  \,\,\, \arc \in \Inedges_\node,
                                                                                              &&
                                                                                                 \press_\arc(0, t)
                                                                                                 =  \press_\node(t),
                                                                                                 \,\,\,  \arc \in \Outedges_\node.
  \end{align}
\end{subequations}
Here, $q_a$ and $\hat q_a$ denote the mass flow and the volumetric
flow in pipe $\arc$, respectively. They scale with the mass density,
\ie, $q_a= \rho_a v_a \varsigma_a$ and $\hat q_a=q_a/\rho_a$, where
$\varsigma_a=d^2_a \pi/4$ is the cross-sectional area of the pipe. In
case of incompressibility, it holds that $\hat{q}_a(x,t)=\hat{q}_a(t)$
is constant along the pipe. The functions $e_\node$ and $\press_\node$ are
auxiliary variables describing internal energy density and pressure at
node $\node$. Note that the second condition in
\eqref{eq:energy-balance}, namely that the out-flowing energy
densities are identical in all (flow-specific outgoing) pipes, rests upon the
assumption of instant mixing of the in-flowing energy densities.

\subsection{Households, depot, and operational constraints}
\label{sec:network-modeling:househelds-depot}

The network modeling is closed by models for the consumers
(households) and the depot of the network operator.
Quantities associated to the arc~$a$ at
node~$n$ are indicated by the subscript~${a:n}$.

For the consumer at $\arc = (\otherNode, \node) \in \Ac$, where the
nodes $m$ and $n$ belong to the foreflow and backflow part of the
network, respectively (\cf Fig.~\ref{fig:sample-network}), the
following conditions are posed for $t\in \timeint$,
\begin{subequations}
  \label{eq:consumer}
  \begin{align}
    \Power_\arc(t)
    & = \volflo_\arc(t) \Delta \edens_{\arc}(t),
    & \vel_\arc(t)
    & \geq 0,
    & \Delta \edens_{\arc}(t)
    & = \edens_{\arc:\otherNode}(t) - \edens_{\arc:\node}(t),
      \label{eq:consumer-power-usage}\\
    \temp_{\arc:\node}(t)
    & = \tempbf,
    & \temp_{\arc:\otherNode}(t)
    & \in [\mintempff,\maxtempff],
    & \temp_{\arc:\otherNode}(t) - \temp_{\arc:\node}(t)
    & \leq \tempdiffhs^\text{c},
    \label{eq:consumer-temperature}\\
    \press_{\arc:\node}(t)
    & \in [\minpresbf,\maxpresbf],
    & \press_{\arc:\otherNode}(t)
    & \in [\minpresff,\maxpresff],
    & \press_{\arc:\otherNode}(t) -\press_{\arc:\node}(t)
    & \in [\Delta p^{\text{c}}_-, \Delta p^{\text{c}}_+].
      \label{eq:consumer-pressure}
  \end{align}
\end{subequations}
The prescribed power consumption~$\Power_\arc$ of the household is
realized by the product of the energy density difference at the arc
and the volumetric flow in~\eqref{eq:consumer-power-usage}. Moreover, the
underlying flow velocity has a pre-specified direction. The consumer's
outflow temperature is set to be equal to the contractually agreed
temperature $\tempbf$. Moreover, the operational constraints ensure a
certain temperature range at each consumption point and define a
maximal temperature difference between foreflow and backflow part of
the consumers.  In addition, minimal and maximal values for the
pressure level at both backflow and foreflow part of the consumer arcs
are prescribed. Finally, the pressure difference between foreflow and
backflow part is bounded.

The depot~$\ad=(\otherNode,\node)$ for operating the district heating
network is modeled by the following conditions for $t \in \timeint$:
\begin{subequations}
  \label{eq:distr-heat-depot}
  \begin{align}
    \label{eq:distr-heat-depot-flow-energy}
    \edens_{\ad:\node}(t) &= u^\text{e}(t), \quad
                            \temp_{\ad:\node}(t) \leq T^\text{net},  \hspace*{2.2cm}
                            \vel_{\ad}(t)  \geq 0, \\
    \label{eq:distr-heat-depot-pressure}
    \press_{\ad:\otherNode}(t) & = u^\text{p}(t), \quad
                                 \press_{\ad:\node}(t)  = \press_{\ad:\otherNode}(t) + u^{\Delta\text{p}}(t).
  \end{align}
\end{subequations}
Here, $u^\text{p}$ prescribes the so-called stagnation pressure of the
network and $u^{\Delta\text{p}}$ is the realized pressure increase at
the depot. The energy density injected at the depot to the foreflow
part of the network is denoted by $u^\text{e}$. The resulting
temperature is bounded above by $T^\text{net}$, which also acts as
temperature limit for all network nodes.

In addition to the operational constraints in \eqref{eq:consumer} and
\eqref{eq:distr-heat-depot}, the pressure in all network nodes is
bounded, \ie, $\press_{\node}(t) \leq p^\text{net}$ for $\node \in
\Vertices$ and $t \in \timeint$.


\section{Port-Hamiltonian formulation of a semi-discrete network
  model}
\label{sec:ph_semidiscrete}

In this section we present a spatially semi-discrete model variant for the
district heating network and discuss its formulation in the
port-Hamiltonian context. Making use of the different hydrodynamic and
thermal time scales, a finite volume upwind discretization yields a
port-Hamiltonian descriptor system for the internal energy density,
in which the solenoidal flow field acts as a time-varying parameter.

We describe the network by means of the following partial
differential-algebraic system for $t\in \timeint$,
\begin{subequations}
  \label{eq:e_net}
  \begin{align}
    \partial_t \edens_\arc
    &= -v_\arc \partial_x \edens_\arc,
    && \arc \in \pipes,
       \label{eq:e_upw_advec}\\
    \edens_\arc(0,t)
    &= \edens_\node(t),
    && \arc \in \outflowArcs_\node(t),
       \qquad \sum_{\arc \in \Inedges_\node} \volflo_\arc
       \edens_\arc(\length_\edge, t) = \sum_{\arc \in \Outedges_\node}
       \volflo_\arc \edens_\arc(0,t),
       \quad \node \in \Vertices,
       \label{eq_mix}\\
    \edens_{\arc:\node}(t)
    &= \edens^\text{bf},
    && \arc \in \Ac,
       \label{eq:e_bc_II}\\
    \edens_{a:\node}(t)
    &= u^\text{e}(t),
    && \arc = \ad,
       \label{eq:e_bc_I}\\
    g(e,v,p)
    &=0 \label{eq:e_upw_algeb}.
  \end{align}
\end{subequations}
This system results from the incompressible pipe model
in~\eqref{eq:distr-heat-pipe-smallacc} and neglecting the cooling term in
the energy balance (\ie, $k_\mathrm{w} =0$). Here, the condition on the backflow temperature for the consumers is expressed in terms of the internal energy density,
\cf, $e^\text{bf}=e(T^\text{bf})$ in~\eqref{eq:e_bc_II}. In the formulation
we use the separation of thermal and hydrodynamic effects and state
the temporal advection of the internal energy density with respect to
the algebraic equations covering the hydrodynamics. So, $g(e,v,p)=0$
in \eqref{eq:e_upw_algeb} contains the hydrodynamic pipe equations,
the pressure continuity at the nodes \eqref{eq:pressure-continuity},
the condition on the households' power consumption
\eqref{eq:consumer-power-usage}, the pressure conditions at the depot
\eqref{eq:distr-heat-depot-pressure}, and the conservation of volume
\begin{align}\label{eq:cons_vol}
  \sum_{\otheredge \in \Inedges_\node} \volflo_\otheredge(t) =
  \sum_{\edge \in \Outedges_\node} \volflo_\edge(t), \quad \node\in
  \Vertices.
\end{align}
Considering the volume balance~\eqref{eq:cons_vol} instead of the
mass balance~\eqref{eq:mass-balance} is very convenient in the
incompressible setting, since the velocity field and hence the induced
volumetric flow are constant along a pipe. Moreover, this description
naturally fits the numerical method of finite volumes.

For the spatial discretization of the hyperbolic-like system
\eqref{eq:e_net} we apply a classical finite volume upwind scheme
\cite{leveque_numerical_2008}. Let $\pipind\in \pipes$, $\pipind\in
\outflowArcs_\node(t_0)$, $\node \in \Vertices$, and consider an
equidistant mesh of cell size $\Delta x_\pipind$, then
\begin{align*}
  \frac{\mathrm{d}}{\mathrm{d} t} \edens_{\pipind,\cellind}
  &=
    -\frac{v_\pipind}{\Delta x_\pipind}(\edens_{\pipind,\cellind} -
    \edens_{\pipind,\cellind-1}),\quad
    \cellind \in \cellset_\pipind,
  \\
  \edens_{\pipind,0}
  &=\edens_\node,
    \quad
    \edens_\node= \frac{\sum_{\otheredge \in \inflowArcs_\node }
    \volflo_\otheredge
    \,\edens_{\otheredge,|\cellset_\otheredge|}}{\sum_{\edge \in
    \outflowArcs_\node} \volflo_\edge},
\end{align*}
where $\edens_{\pipind,\cellind}$ denotes the internal energy density
with respect to the finite volume cell~$\cellind$ of pipe~$\pipind$
with cell index set $\cellset_\pipind$. For the first cell
($\cellind=1$) we make use of the quantity at the node that results
from \eqref{eq_mix}. We summarize the unknown energy densities in a
vector $e=(e_1,...,e_\kappa)^T$, $e_{f(\alpha,\beta)}=e_{\alpha,\beta}$ by ordering pipe- and cell-wise according to the mapping
$f(\pipind,\cellind) = \cellind + \sum_{k=1}^{\pipind-1}
|\cellset_k|$, $ \pipind \in \pipes$,  $\cellind \in
\cellset_\pipind$, in particular $\numcells=\sum_{\pipind \in \pipes}
|\cellset_\pipind|$.
Then, a semi-discrete version of the network model~\eqref{eq:e_net} is
given by the following descriptor system
\begin{align}\label{eq:upw:sys}
  &\frac{\mathrm{d}}{\mathrm{d} t} e= A(\vel) \,e + B(\vel)\, u,
    \quad y = C e,
  \\ \nonumber
  &\text{subject to } \vel=G(e).
\end{align}
The system matrices $A(w)\in \mathbb{R}^{\numcells\times \numcells}$
and $B(w)\in \mathbb{R}^{\numcells \times 2}$ can be interpreted as
parameter-dependent quantities, where the (vector-valued) parameter
$w$ represents a spatially discretized solenoidal volume-preserving
velocity field. So,
\begin{equation*}
  A_{f(\pipind,\cellind),f(\othpipind,\othcellind)}(w)= \partial
  \frac{\mathrm d}{\mathrm{d} t} \edens_{\pipind,\cellind}(w)/{\partial
    \edens_{\othpipind,\othcellind}}
\end{equation*}
holds. The special velocity field
belonging to the hydrodynamic network equations~\eqref{eq:e_upw_algeb}
is formally stated as $v=G(e)$. We assume a setting in which $v$ is
time-continuous. In~\eqref{eq:upw:sys} the input~$u$ consists of the
energy densities $u^\text{e}$ injected at the depot into the foreflow
part and $\edens^{\text{bf}}$ returning from the consumers into the
backflow part of the network,
$u=(u^\text{e},\edens^{\text{bf}})^T \in \mathbb{R}^{2} $. The
output~$y$ typically refers to energy densities in pipes supplying the
consumers, implying $C\in \mathbb{R}^{c \times \kappa}$.

\begin{theorem}\label{theo:semi_discr_pH}
  Let $w$ be a (spatially discretized) solenoidal volume-preserving
  time-continuous velocity field. Then, the semi-discrete network
  model~\eqref{eq:upw:sys} can be embedded into a family of
  parameter-dependent port-Hamiltonian systems
  \begin{equation}\label{eq:pH_net}
    \frac{\mathrm{d}}{\mathrm{d}t} e
    = (J(w)-R(w)) Q e + \tilde B(w) \tilde u,
    \quad
    \tilde y = \tilde{B}^T(w) Q e,
  \end{equation}
  with $\tilde u=(u^T,0,\dotsc,0)^T\in \mathbb{R}^{2+c}$ which contains
  the original outputs as subset.
\end{theorem}

\begin{remark}
  Theorem~\ref{theo:semi_discr_pH} implies that there exists an energy
  matrix~$\ematpH$ such that
  \begin{equation}\label{eq:upw_lyap}
    \ematpH A(w) + A^T(w) \ematpH \leq 0
  \end{equation}
  for all solenoidal volume-preserving velocity fields $w$. Thus, the
  Hamiltonian $\mathcal{H}(\edens) = \edens^T \ematpH \edens$ is a
  Lyapunov function for the parameter-dependent system
  \cite{antoul_approx}. The energy matrix $\ematpH$ can be particularly
  constructed as a diagonal matrix with positive entries, \ie, $
  \ematpH_{f(\pipind,\cellind),f(\pipind,\cellind)} = 
  \varsigma_\pipind\, \Delta x_\pipind $ for $\pipind \in \pipes$, $\cellind \in V_\pipind$, where
  $  \varsigma_\pipind\,\Delta x_\pipind$ is the volume of each discretization
  cell in pipe $\pipind$.

  Note that a change of the flow direction, which might occur in case
  of cycles, yields a structural modification of the system matrix
  $A(w)$, but does not affect the
  stability of the system. However, it might cause a discontinuity in
  the velocity field such that \eqref{eq:upw:sys},  or
  \eqref{eq:pH_net} respectively, only allows for a weak solution.
\end{remark}

\begin{proof}[of~Theorem~\ref{theo:semi_discr_pH}]
  Let the positive definite diagonal matrix $Q\in\mathbb{R}^{\kappa
    \times \kappa}$ with
  $\ematpH_{f(\pipind,\cellind),f(\pipind,\cellind)} =
    \varsigma_\pipind\,\Delta x_\pipind  > 0$ be given. Then, we define the
  matrices $J$ and $R$ by
  \begin{equation*}
    J(w) = \frac{1}{2}(A(w)\ematpH^{-1} -(A(w)\ematpH^{-1})^T),
    \quad
    R(w) = -\frac{1}{2}(A(w)\ematpH^{-1} + (A(w)\ematpH^{-1})^T).
  \end{equation*}
  Obviously, $A(w)=(J(w)-R(w))Q$ holds. The properties $J=-J^T$ and
  $R=R^T$ of port-Hamiltonian system matrices are satisfied by
  construction for any parameter $w$. The positive semi-definiteness of~$R$ follows from
  the Lyapunov inequality~\eqref{eq:upw_lyap}.
  Considering
  \begin{equation*}
    L(w) = \ematpH A(w) + A^T(w) \ematpH,
    \quad L_{f(\pipind,\cellind),f(\pipind,\cellind)}(w) = -2
    \ematpH_{f(\pipind,\cellind),f(\pipind,\cellind)}
    \frac{w_\pipind}{\Delta x_\pipind} = -2\volflo_\pipind\leq 0,
  \end{equation*}
  the volume-preservation of $w$ ensures that the symmetric matrix
  $L(w)$ is weakly diagonal dominant. Hence, $L(w)$ is negative
  semi-definite, yielding
  \begin{equation*}
    x^T R(w) x= -\frac{1}{2} (Q^{-1}x)^T \, L(w) \, (Q^{-1}x) \geq 0
    \quad\text{for all}\quad
    x \in \mathbb{R}^\kappa.
  \end{equation*}
  Here, $R(w)$ acts as the passivity matrix since the system has no feed-through term. The port matrix
  $\tilde{B}(w)\in \mathbb{R}^{\kappa \times 2+c}$ defined by
  \begin{align*}
    \tilde{B}(w)=[B(w), \,\, (CQ^{-1})^T]
  \end{align*}
  ensures that the outputs of the network model are contained in the
  output set of the port-Hamiltonian system, \ie, $\tilde{B}^T(w)Q=
  [B^T(w)Q, \,\,C]^T$. Finally note that the parameter-dependent
  port-Hamiltonian system matrices $J(w)$, $R(w)$, and $\tilde B(w)$
  are continuous in time due to the given time-regularity of the
  parameter $w$.
\end{proof}

\begin{remark}
  We point out that applying the stated framework to the other pipe
  models presented in Sect.~\ref{subsec:pipe-modeling} is
  non-trivial. Already the consideration of the cooling term in the
  energy balance, \cf pipe model \eqref{eq:distr-heat-pipe-smallacc},
  which acts dissipative requires a generalization of the
  port-Hamiltonian description. We refer to Sect.~\ref{sec:pH-modeling}
  for an infinite-dimensional port-Hamiltonian formulation of the
  compressible thermodynamic pipe flow \eqref{eq:distr-heat-pipe}.
\end{remark}


\section{Numerical study on network operation}
\label{sec:case_study}

In this section we demonstrate the potential for optimization of
district heating networks. Operating the network according to certain
exogenously given temporal profiles for the internal energy densities
injected at the depot may lead to high amplitudes in the feed-in
power. The avoidance of such power peaks in the feed-in prevents that
using additional energy sources, such as gas storages, is required for
covering the heating demand of the consumers. This is environmental
friendly, while saving resources and operational costs.

\begin{figure}[t]
  \centering
  \includegraphics[scale=0.7]{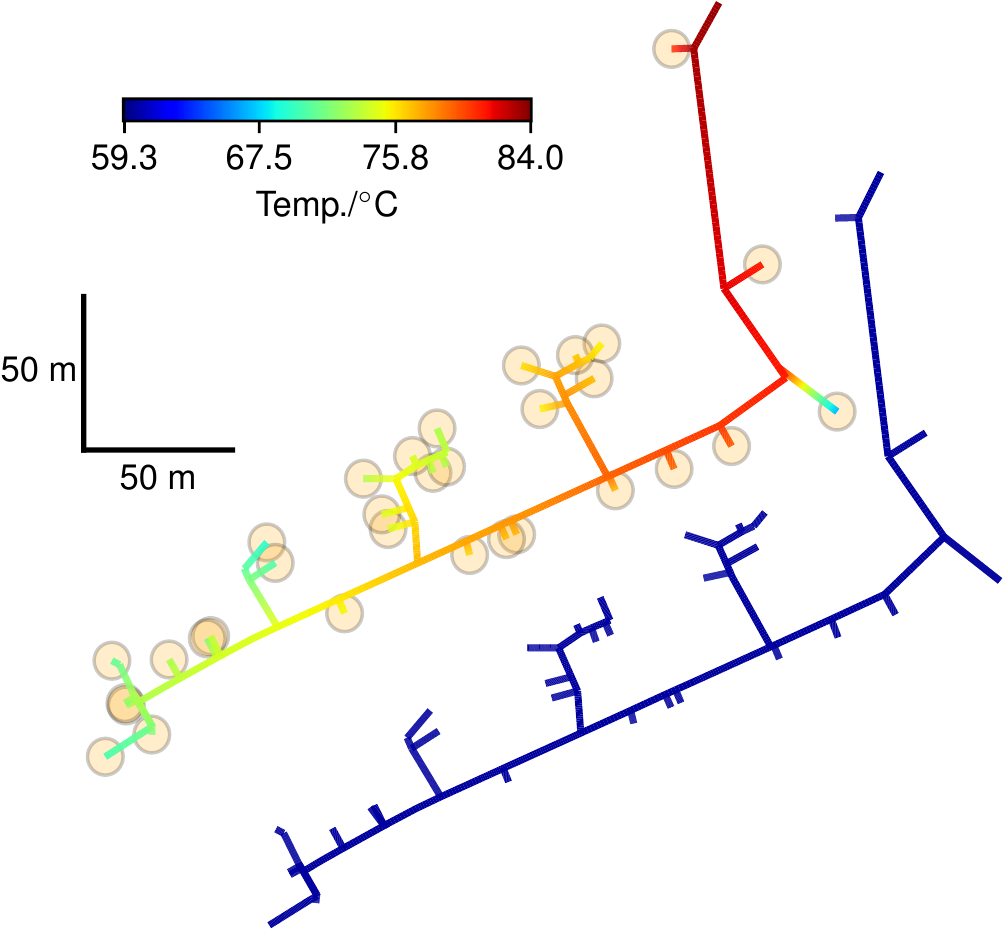}
  \caption{Real-world heating network supplying several streets. The
    network consists of the foreflow part (top) and the backflow part
    (bottom), where the households are indicated
    by circles. The topology has been provided by Technische Werke
    Ludwigshafen AG, Germany. The color plot visualizes a simulated
    temperature distribution for a certain time~$t^\star$, where
    $T(u^\text{e}(t^\star))=\SI{84}{\degreeCelsius}$. The backflow
    temperature is constant at $T^\text{bf}=\SI{60}{\degreeCelsius}$ due to the
    use of the network model \eqref{eq:e_net} where cooling effects are neglected.}
  \label{fig:topol_street}
\end{figure}

In the numerical case study we employ a real-world district heating
network supplying different streets by means of the port-Hamiltonian
semi-discrete network model~\eqref{eq:upw:sys}. For the time
integration we use an implicit midpoint rule with constant time
step~$\Delta t$. The topology of the network and the data of the
pipelines come from the Technische Werke Ludwigshafen AG; see
Fig.~\ref{fig:topol_street} and Table~\ref{tab:runtime_fullred}.
For the presented simulation, a time
horizon of \SI{50}{\hour} is studied. The consumption behavior of the
households is modeled by standardized profiles used in the operation
of district heating networks \cite{SLP_BGW}
for a mean environmental temperature of~\SI{3}{\degreeCelsius}. The total consumption of
all households is \SI{108}{\kilo\watt} on temporal average and rises up to a
maximum of~\SI{160}{\kilo\watt}. Given the internal energy density~$u^\text{e}$
injected at the depot as input, the feed-in power can be considered as
the response of the network system, \ie,
\begin{equation*}
  \feedin = (u^\text{e}-e_{\ad:m} )\sum_{\arc \in \Ac} \volflo_\arc.
\end{equation*}
Note that due to the neglect of cooling in~\eqref{eq:e_net},
$e_{\ad:m}=e(T^\text{bf})$ holds, where the backflow temperature at
the consumers is fixed here to $T^\text{bf}=\SI{60}{\degreeCelsius}$.

\begin{table}[t]
  \centering
  \caption{Graph-associated outline data for the street network in
    Fig.~\ref{fig:topol_street}. The total pipe length of the foreflow
    part is \SI{835.5}{\meter} and of the backflow part \SI{837.0}{\meter}.}
  \label{tab:runtime_fullred}
  \begin{tabular}{cccccc}
    \toprule
    Pipes $|\pipes|$
    & Consumers $|\Ac|$
    & Depot
    & Arcs $|\Arcs|$
    & Nodes $|\Vertices|$ & Loops
    \\
    \midrule
    162 & 32 & 1 & 195 & 162 & 2 \\
    \bottomrule
  \end{tabular}
\end{table}

\begin{figure}[t]
  \centering
  \includegraphics[scale=0.7]{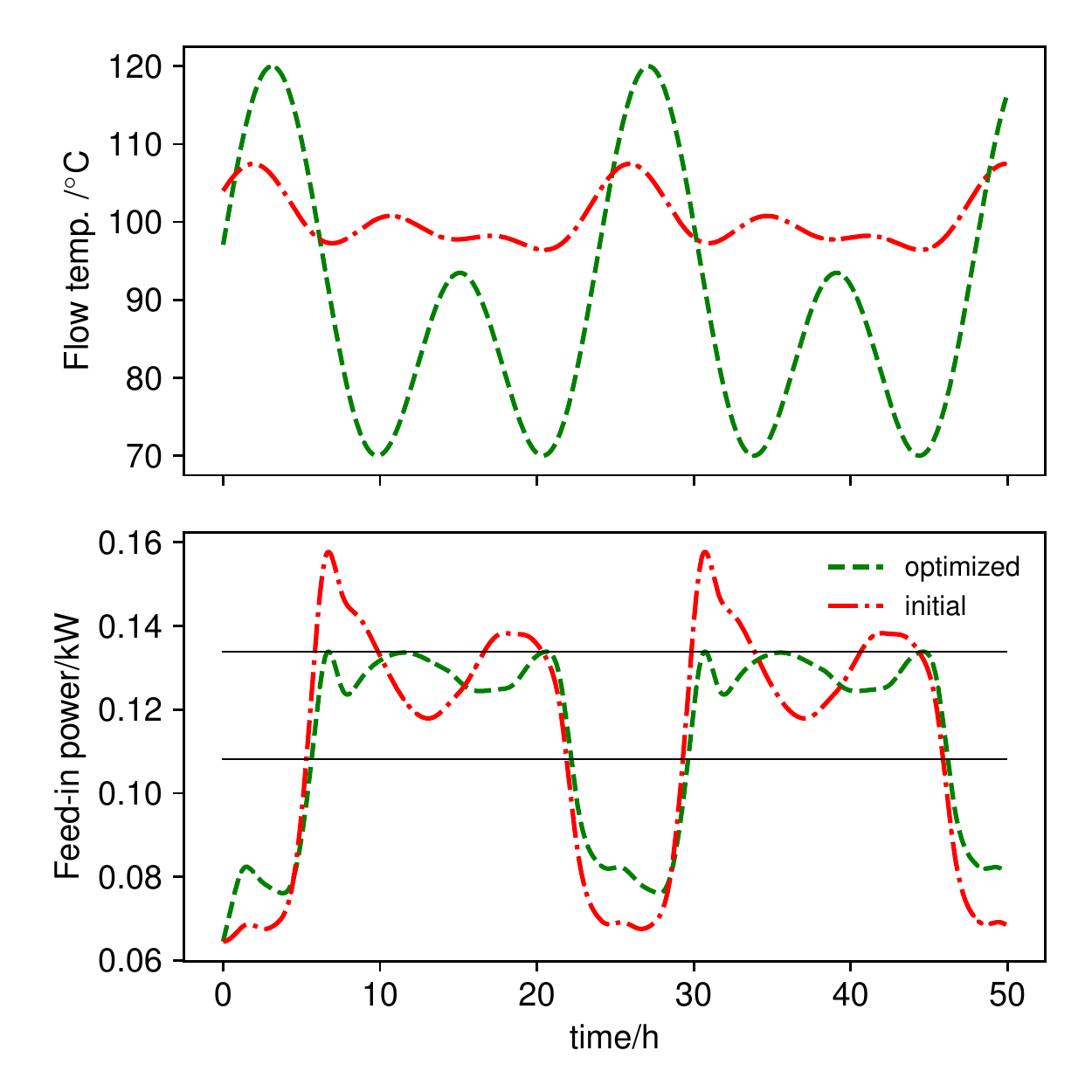}
  \caption{Flow temperature at depot $T(u^\text{e})$ (top) and
    corresponding feed-in power (bottom) over time for two different
    injection profiles marked in dashed-dotted red and dashed green,
    $\Delta t= \SI{5}{\minute}$. The upper solid, black line indicates the power
    threshold $\Plim$, the lower one the mean feed-in power over time.}
  \label{fig:opt_street}
\end{figure}

The traveling time of the heated water from the depot to the consumers
(households) allows to choose from different injection profiles, when
covering the aggregated heating demand in the
network. Figure~\ref{fig:opt_street} shows the injected temperature
$T(u^\text{e})$ and the corresponding feed-in power for two different
input profiles. Supplying an almost constant energy density~$u^\text{e}$ over time
yields pronounced power peaks (dashed-dotted red curves). These
undesired peaks can be avoided when using an input that is varying in
time with respect to the expected consumer demands. For the
illustrated improved input, the feed-in power is bounded by
$\bar{P}_\text{in} = \SI{134}{\kilo\watt}$ (dashed green curves). This promising result asks for a rigorous optimal control of the network in further studies.

\section{Port-Hamiltonian formulation of compressible thermodynamic pipe flow}
\label{sec:pH-modeling}

The adequate handling of thermal effects requires the generalization of the port-Hamiltonian framework by combing the Hamiltonian with an entropy function. In this section we embed the partial differential model \eqref{eq:distr-heat-pipe} for a compressible thermodynamic turbulent pipe flow into the GENERIC-formalism, which has lately been studied in \cite{BadMBM18,BadZ18}, and present an infinite-dimensional thermodynamically consistent port-Hamiltonian description.

The thermodynamic pipe flow model  \eqref{eq:distr-heat-pipe} can be reformulated as a generalized (non-linear) port-Hamiltonian system in operator form for $z=(\rho, M, e)^T$, $M=(\rho v)$,
\begin{equation}
\label{eq:GENERIC_ph}
\begin{aligned}
\frac{\mathrm{d}z}{\mathrm{d}t}&=\left(\mathcal{J}(z)-\mathcal{R}(z)\right)\frac{\delta\mathcal{E}(z)}{\delta z}+\mathcal{B}(z)u(z)\quad &\text{in }  \mathcal{D}_z^*,\\
y(z)&=\mathcal{B}^*(z)\frac{\delta\mathcal{E}(z)}{\delta z} \quad &\text{in } \mathcal{D}_u^*,
\end{aligned}
\end{equation}
where  $\mathcal{Z}=\{z\in \mathcal{D}_z \,|\,\rho\geq\delta \text{ with } \delta>0 \text{ almost everywhere} \}\subset \mathcal{D}_z$ denotes the state space with the Sobolev space  $\mathcal{D}_z=W^{1,3}((0,\ell);\mathbb{R}^3)$ being a reflexiv Banach space. For $z\in \mathcal{Z}$ the operators $\mathcal{J}(z)[\cdot]$, $\mathcal{R}(z)[\cdot]:\mathcal{D}_z \rightarrow \mathcal{D}_z^*$ are linear and continuous, moreover $\mathcal{J}(z)$ is skew-adjoint and $\mathcal{R}(z)$ is self-adjoint semi-elliptic, i.e., $\langle\varphi,\mathcal{J}(z)\psi\rangle=-\langle\psi,\mathcal{J}(z)\varphi\rangle$ and $\langle\varphi,\mathcal{R}(z)\psi\rangle=\langle\psi,\mathcal{R}(z)\varphi\rangle\geq 0$ for all $\varphi,\psi \in\mathcal{D}_z$. The system theoretic input is given by $u(z)\in \mathcal{D}_u= L^{q}(\{0,\ell\})$ with linear continuous operator $\mathcal{B}(z)[\cdot]:\mathcal{D}_u\rightarrow\mathcal{D}_z^*$ and dual space $\mathcal{D}_u^*= L^{p}(\{0,\ell\})$, ${1}/{q}+{1}/{p}=1$. The system theoretic output is denoted by $y(z)$.  The form of the energy functional $\mathcal{E}$ and the port-Hamiltonian operators $\mathcal{J}(z)$, $\mathcal{R}(z)$ and $\mathcal{B}(z)[\cdot]$ are derived as follows.

\begin{remark}
We assume that all relevant mathematical statements hold for an arbitrary but fixed time parameter $t\in (t_0,t_\text{end}]$. The function spaces $\mathcal{D}_z$ and $\mathcal{D}_u$ associated with the spatial evolution are chosen in an ad-hoc manner, \ie, we assume that the considered fields and functions satisfy certain regularity requirements.  A mathematically rigorous justification requires an analytical consideration of the generalized port-Hamiltonian system. The corresponding functional analytical and structural questions are the focus of ongoing work.
\end{remark}

Accounting for the thermodynamic behavior of the pipe flow, \eqref{eq:GENERIC_ph} is composed of a Hamiltonian and a generalized gradient system. This is reflected in the energy functional that is an exergy-like functional consisting of a Hamiltonian and an entropy part, i.e.,
\begin{align*}
\mathcal{E}(z)=\mathcal{H}(z)-\vartheta\mathcal{S}(z), \quad
\mathcal{H}(z)=\int_0^\ell \left(\frac{|M|^2}{2\rho}+e+\rho gh\right)\mathrm{d}x, \quad \mathcal{S}(z)=\int_0^\ell s(\rho,e)\,\mathrm{d}x.
\end{align*}
where the outer ground temperature $\vartheta$ is assumed to be constant.
Introducing the ballistic free energy $H(\rho,e)= e-\vartheta s(\rho,e)$ \cite{Feireisl12}, the functional $\mathcal{E}$  and its variational derivatives become 
\begin{align*}
\mathcal{E}(z)&=\int_0^\ell \left(\frac{|M|^2}{2\rho}+H(\rho,e)+\rho gh\right)\mathrm{d}x\\
\frac{\delta\mathcal{E}(z)}{\delta z}&=
\left( \frac{\delta\mathcal{E}(z)}{\delta \rho}, \frac{\delta\mathcal{E}(z)}{\delta M},\frac{\delta\mathcal{E}(z)}{\delta e} \right)^T
=\left(\left(-\frac{|M|^2}{2\rho^2}+\frac{\partial H}{\partial\rho}+gh\right), \frac{M}{\rho}, \frac{\partial H}{\partial e}\right)^T.
\end{align*}
The port-Hamiltonian operators in \eqref{eq:GENERIC_ph} are assembled with respect to the (block-) structure of the state $z$. Let $\varphi, \psi \in \mathcal{D}_z$ be two block-structured test functions, i.e., $\varphi = (\varphi_\rho, \varphi_M, \varphi_e)^T$. Then the skew-adjoint operator $\mathcal{J}(z)$ is given by
\begin{subequations}
\label{eq:Skew_Op_J}
\begin{align}
\mathcal{J}(z)=\left[ \begin{array}{ccc}0 & \mathcal{J}_{\rho,M}(z) & 0 \\ \mathcal{J}_{M,\rho}(z) & \mathcal{J}_{M,M}(z) &\mathcal{J}_{M,e}(z) \\ 0 & \mathcal{J}_{e,M}(z) & 0\end{array}\right],
\end{align}
associated with the bilinear form
\begin{align*}
\langle\varphi,\mathcal{J}(z)\psi\rangle & =
\langle\varphi_\rho,\mathcal{J}_{\rho,M}(z)\psi_M\rangle+
\langle\varphi_M,\mathcal{J}_{M,\rho}(z)\psi_\rho\rangle+
\langle\varphi_M,\mathcal{J}_{M,M}(z)\psi_M\rangle\\ &\quad +
\langle\varphi_M,\mathcal{J}_{M,e}(z)\psi_e\rangle+
\langle\varphi_e,\mathcal{J}_{e,M}(z)\psi_M\rangle.
\end{align*}
Its entries are particularly defined by the following relations,
\begin{align}
\langle\varphi_\rho,\mathcal{J}_{\rho,M}(z)\psi_M\rangle&=-\langle\psi_M,\mathcal{J}_{M,\rho}(z)\varphi_\rho\rangle=\int_0^\ell\rho(\psi_M\partial_x)\varphi_\rho\,\mathrm{d}x,\\
\langle\varphi_M,\mathcal{J}_{M,M}(z)\psi_M\rangle&=-\langle\psi_M,\mathcal{J}_{M,M}(z)\varphi_M\rangle=\int_0^\ell M((\psi_M\partial_x)\varphi_M-(\varphi_M\partial_x)\psi_M)\,\mathrm{d}x,\\
\langle\varphi_e,\mathcal{J}_{e,M}(z)\psi_M\rangle&=-
\langle\psi_M,\mathcal{J}_{M,e}(z)\varphi_e\rangle=\int_0^\ell e(\psi_M\partial_x)\varphi_e+(\psi_M\partial_x)(\varphi_ep)\,\mathrm{d}x \label{eq:ph-J}
\end{align}
\end{subequations}
that result from the partial derivatives in \eqref{eq:distr-heat-pipe}. The self-adjoint semi-elliptic operator $\mathcal{R}(z)$ is composed of two operators that correspond to the friction in the pipe $\mathcal{R}^\lambda(z)$ and the temperature loss through the pipe walls $\mathcal{R}^{k_\mathrm{w}}(z)$. It is given by
\begin{subequations}
\label{eq:Dis_Op_R}
\begin{align}
\mathcal{R}(z)=\mathcal{R}^\lambda(z)+\mathcal{R}^{k_\mathrm{w}}(z)=\left[\begin{array}{ccc}0 & 0 & 0 \\ 0 & \mathcal{R}^\lambda_{M,M}(z) &\mathcal{R}^\lambda_{M,e}(z) \\ 0 & \mathcal{R}^\lambda_{e,M}(z) & \mathcal{R}^\lambda_{e,e}(z)+\mathcal{R}^{k_\mathrm{w}}_{e,e}(z)\end{array}\right],
\end{align}
associated with the bilinear form,
\begin{align*}
\langle\varphi,\mathcal{R}(z)\psi\rangle & =
\langle\varphi_M,\mathcal{R}^\lambda_{M,M}(z)\psi_M\rangle+
\langle\varphi_M,\mathcal{R}^\lambda_{M,e}(z)\psi_e\rangle+
\langle\varphi_e,\mathcal{R}^\lambda_{e,M}(z)\psi_M\rangle\\ &\quad+\langle\varphi_e,(\mathcal{R}^\lambda_{e,e}(z)+\mathcal{R}^{k_\mathrm{w}}_{e,e}(z))\psi_e\rangle.
\end{align*}
Its entries are
\begin{align}\label{eq:ph-R1}
\langle\varphi_M,\mathcal{R}^\lambda_{M,M}(z)\psi_M\rangle&=\int_0^\ell \varphi_M\left(\frac{\lambda}{2d}\frac{\temp}{\vartheta}\rho|\vel|\right)\psi_M\,\mathrm{d}x,\\
\langle\varphi_M,\mathcal{R}^\lambda_{M,e}(z)\psi_e\rangle=
\langle\psi_e,\mathcal{R}^\lambda_{e,M}(z)\varphi_M\rangle&=\int_0^\ell -\varphi_M\left(\frac{\lambda}{2d}\frac{\temp}{\vartheta}\rho|\vel|\vel\right)\psi_e\,\mathrm{d}x,\\
\langle\varphi_e,(\mathcal{R}^\lambda_{e,e}(z)+\mathcal{R}^{k_\mathrm{w}}_{e,e}(z))\psi_e\rangle&=\int_0^\ell \varphi_e\left(\frac{\lambda}{2d}\frac{\temp}{\vartheta}\rho|\vel|\vel^2+\frac{4k_\mathrm{w}}{d}\temp\right)\psi_e\,\mathrm{d}x.\label{eq:ph-R3}
\end{align}
\end{subequations}
Note that the state dependencies of pressure $p=p(\rho,e)$ and temperature $T=T(\rho,e)$ occurring in \eqref{eq:ph-J} and \eqref{eq:ph-R1}-\eqref{eq:ph-R3} are prescribed by the state equations, cf.\ Remark~\ref{rem:state}. Moreover, $v=M/\rho$ and $\lambda=\lambda(v,e)$ hold for the velocity and the friction factor, respectively. Assuming consistent state equations, \eg, ideal gas law, \cf Remark~\ref{rem:ideal-gas}, the operators in \eqref{eq:Skew_Op_J} and \eqref{eq:Dis_Op_R} fulfill the non-interacting conditions
\begin{align*}
\mathcal{J}(z)\frac{\delta\mathcal{S}(z)}{\delta z}=0, \quad \quad \mathcal{R}^\lambda(z)\frac{\delta\mathcal{H}(z)}{\delta z}=0,
\end{align*}
which arise in the GENERIC context \cite{BadMBM18,BadZ18} and ensure that the flows of the Hamiltonian and the gradient system do not overlap. Finally, concerning the system theoretic input and output, the state dependent input is given as $u(z)\in \mathcal{D}_u$ by $u(z)=[M/\rho]|^\ell_0$. Then, the port operator $\mathcal{B}(z)[\cdot]:\mathcal{D}_u\rightarrow \mathcal{D}_z^\star$ is specified through the pairing 
\begin{align*}
\langle\varphi,\mathcal{B}(z)u(z)\rangle=-\left.\left[(\varphi_\rho\rho+\varphi_MM+\varphi_e(e+p))\,u(z)\right]\right|_0^\ell,
\end{align*}
which originates from the boundary terms, when applying partial integration to parts of (\ref{eq:distr-heat-pipe}).
With the adjoint operator $\mathcal{B}^*(z)[\cdot]:\mathcal{D}_z\rightarrow \mathcal{D}_u^*$, i.e., 
$\langle\varphi,\mathcal{B}(z)u(z)\rangle=\langle\mathcal{B}^*(z)\varphi,u(z)\rangle $, the system theoretic output reads
\begin{align*}
y(z)=\mathcal{B}^*(z)\frac{\delta\mathcal{E}(z)}{\partial z}=-\left.\left[\frac{|M|^2}{2\rho}+p+H(\rho,e)+\rho gh\right]\right|^\ell_0.
\end{align*}

\begin{remark}\label{rem:ideal-gas}
In the port-Hamiltonian framework the choice of the state variables in the interplay with the energy functional is crucial for encoding the physical properties in the system operators. Hence, asymptotic simplifications as, \eg, the limit to incompressibility in the hydrodynamics \eqref{eq:distr-heat-pipe-incomp}, are not straightforward, since they change the underlying equation structure. However, system~\eqref{eq:GENERIC_ph} is well suited when, \eg, dealing with gas networks. Then, it can be closed by using, \eg, the ideal gas law, implying \begin{align*}
s(\rho,e)=\frac{R}{2} \rho \operatorname{ln}\left(\heatcap\frac{e^3}{\rho^5}\right), \quad \temp(\rho,e)=\frac{2}{3R}\frac{e}{\rho}, \quad  \press(\rho,e)=\frac{2}{3}e,
\end{align*}
with specific gas constant $R$ and heat capacity $\heatcap$.
\end{remark}

\section{Research perspectives}
\label{sec:concl-persp}

An energy-based port-Hamiltonian framework is very suitable for
optimization and control when dealing with subsystems coming from
various different physical domains, such as hydraulic, electrical, or
mechanical ones, as it occurs when coupling a district heating
network with a power grid, a waste incineration plant, or a gas
turbine. The formulation is advantageous as it brings different scales
on a single level, the port-Hamiltonian character is inherited by the
coupling, and the physical properties are directly encoded in the
structure of the equations. However, to come up with efficient
adaptive optimization strategies based on port-Hamiltonian model
hierarchies for complex application issues on district heating
networks, there are still many mathematical challenges to be handled.

In this paper we contributed with an infinite-dimensional
and thermodynamically consistent formulation for a compressible turbulent pipe flow,
which required to set up a (reversible) Hamiltonian system and a
generalized (dissipative) gradient system with suitable degeneracy
conditions. In particular, the choice of an appropriate energy
function was demanding. The asymptotic transition to an incompressible
pipe flow is non-trivial in this framework, since it changes the
differential-algebraic structure of the equations and hence requires
the reconsideration of the variables and the modification of the
energy function. In view of structure-preserving discretization and
model reduction the use of Galerkin projection-based techniques seems
to be promising. However, the choice of the variables and the formulation
of the system matrices crucially determine the complexity of the
numerics as, \eg, the works \cite{CBG2016,Egger:2018,c47:liljegren-sailer:2019} show. Especially, the handling of the nonlinearities requires adequate complexity-reduction strategies. Interesting to explore are certainly also structure-preserving time-integration schemes, see, \eg, \cite{KotL18,MorM19}.

In the special case of the presented semi-discrete district heating network model that makes use of the different hydrodynamic and thermal time scales and a suitable finite volume upwind discretization we came up with a finite-dimensional port-Hamiltonian system for the internal energy density where the solenoidal flow field acts a time-varying parameter. This system is employed for model reduction (moment matching) in \cite{rein:2019a} and for optimal control in \cite{rein:2019b}.

The application of the port-Hamiltonian modeling framework
for coupled systems leads to many promising ideas for the
optimization of these systems.
Due to the complexity and size of the respective optimization models,
a subsystem-specific port-Hamiltonian modeling together with suitable
model reduction techniques allows for setting up a coupled model
hierarchy for optimization, which paves the way for highly efficient
adaptive optimization methods; \cf, \eg, \cite{MSS18}, where a
related approach has shown to be useful for the related field of
gas network optimization.


\begin{acknowledgements}
  The authors acknowledge the support by the German BMBF, Project
  \emph{EiFer -- Energy efficiency via intelligent heating networks} and are very grateful for the  provision of the data by their industrial partner Technische Werke Ludwigshafen AG. 
  Moreover, the support of the DFG within the CRC~TRR~154, subprojects~A05, B03, and B08, as well as within the RTG~2126 \emph{Algorithmic Optimization} is acknowledged.
\end{acknowledgements}


\bibliographystyle{spmpsci}
\bibliography{mso-dhn}

\end{document}